\documentclass[final]{elsarticle}

\usepackage{amssymb,latexsym}
\usepackage{amsmath}

\biboptions{sort&compress}

\journal{: \; DMTCS}

\begin{document}

\newtheorem{lema}{Lemma}
\newtheorem{teo}{Theorem}
\newproof{prova}{Proof}
\newtheorem{corol}{Corollary}

\begin{frontmatter}

\title{Rapidly converging Ramanujan-type series for Catalan's constant}

\author{F.~M.~S.~Lima}

\address{Institute of Physics, University of Bras\'{i}lia, P.O.~Box 04455, 70919-970, Bras\'{i}lia-DF, Brazil}


\ead{fabio@fis.unb.br}


\begin{abstract}
In this note, by making use of a known hypergeometric series identity, I prove two Ramanujan-type series for the Catalan's constant.  The convergence rate of these central binomial series surpasses those of all known similar series, including a classical formula by Ramanujan and a recent formula by Lupas.  Interestingly, this suggests that an Ap\'{e}ry-like irrationality proof could be found for this constant.
\end{abstract}


\end{frontmatter}

\section{Introduction}

Catalan's constant, so named in honor to E. C. Catalan (1814-1894), who first developed series and definite integrals representations for it, is a classical mathematical constant defined as~\cite{Catalan1865,Finch}
\begin{equation}
G := \sum_{n=0}^\infty{\frac{(-1)^n}{(2n+1)^2}} = 0.9159655941\ldots
\label{eq:def}
\end{equation}
This constant is a special value of some important functions such as the Dirichlet's beta function $\,\beta{(s)} := \sum_{n=0}^\infty{(-1)^n/(2n+1)^s}$, 
 namely $\,\beta{(2)} = G$, and the Clausen's function $\,\mathrm{Cl}_2(\theta) := \Im{\left(\mathrm{Li}_2(e^{i\theta})\right)} = -\int_0^{\theta}{\ln{\!\left[2 \, \sin{(t/2)}\right]} \, dt}$, namely\footnote{As usual, $\mathrm{Li}_2(z)\,$ denotes the dilogarithm function, defined as $\,\sum_{n=1}^\infty{z^n/n^2}$ for real values of $z$, $|z|<1$, and extended to $\mathbb{C}$, except for the cut $[1,\infty)$, by analytic continuation.}
\begin{equation}
\mathrm{Cl}_2\!\left(\frac{\pi}{2}\right) = G \quad \mathrm{and} \quad \mathrm{Cl}_2\!\left(\frac{3 \pi}{2}\right) = -\,G \, .
\label{eq:Clausen}
\end{equation}
It also appears on special values of some other important functions, as e.g.
\begin{eqnarray*}
\psi_1\!\left(\frac14\right) &=& \pi^2 + 8 G \, , \quad \psi_1\!\left(\frac34\right) = \pi^2 - 8 G \, , \quad \mathrm{and} \nonumber \\
\zeta'\!\left(-1,\frac14\right) &=& \frac18 \, \ln{A}  +\frac{G}{4 \pi} -\frac{1}{96} \, ,
\end{eqnarray*}
where $\psi_1(x):= \dfrac{d \psi(x)}{d x}$ is the trigamma function and $\,\psi(x) := \dfrac{d }{d x} \ln{\Gamma(x)} \,$ is the digamma function, $\Gamma(x)\,$ being the classical Euler gamma function. As usual, $\zeta'(s,a)$ is the partial derivative of the Hurwitz zeta function with respect to $s$. $\,A$ is the Glaisher-Kinkelin constant, which in turn is related to the Riemann zeta function by $\,\ln{A} = 1/12 \, -\zeta'{(-1)}$.

For positive integers $\,n$, we can trace an analogy between $\,\beta{(n)}\,$ and $\,\zeta{(n)} := \sum_{k=1}^\infty{1/k^n}$, the Riemann zeta function, since both $\,\zeta(2n)\,$ and $\,\beta(2n-1)\,$ are rational multiples of $\,\pi^{2n}$ and $\,\pi^{2n-1}$, respectively, whereas closed-form expressions for $\,\zeta(2n+1)$ and $\beta(2n)$ in terms of other basic constants are unknown~\cite{Lima2012}.  However, the proof by Ap\'{e}ry (1978) that $\zeta{(3)}$ is irrational~\cite{Poorten} has created an `asymmetry' in that analogy because the irrationality of $\,\beta(2)$, though very suspected, remains unproven.\footnote{It is not known if $\,G\,$ is irrational. Presently, the only known irrationality results for $\,\beta(2n)\,$ are the proofs by Rivoal and Zudilin (2003) that it is irrational for infinitely many values of $\,n$, and that at least one of the seven numbers $\beta(2), \, \ldots , \beta(14)$ is irrational~\cite{RZ}. \smallskip}

From the point of view of numerical computation, Catalan himself (1865) computed $G$ to $14$ decimal places~\cite{Catalan1865}.  By making use of a technique from Kummer, Bresse (1867) computed it to $24$ decimals, a result that was improved to $32$ decimals by Glaisher (1913)~\cite{Wolfram}.  With the advent of computers, $G$ has been computed to a large number of digits. For instance, Yee and Chan (2009) computed it to $31$ billion decimals~\cite{Yee}.  Their computation employs two formulas, one of which is a central binomial series due to Ramanujan (1915)~\cite{Ramanujan}:\footnote{This can be proved from the fact that $\:G = -\int_0^{\pi/4}{\ln{(\tan{\theta})} \: d\theta} = -\frac32 \, \int_0^{\pi/12}{\ln{(\tan{\theta})} \: d\theta}$, as nicely described in Ref.~\cite{Bradley}.}
\begin{equation}
G = \frac{\pi}{8} \, \ln{(2+\sqrt{3})} + \frac38 \, \sum_{n=0}^{\infty}{\frac{1}{\,(2n+1)^2 \, \binom{2n}{n}}} \, .
\label{eq:Rama}
\end{equation}
On searching for similar rapidly converging series, Lupas (2000) has found the following alternating series~\cite{Lupas}
\begin{equation}
G = -\,\frac{1}{64} \, \sum_{n=1}^{\infty}{(-1)^n \, \frac{\, 2^{8n} \, (40n^2-24n+3)}{\,n^3 \, (2n-1) \, \binom{2n}{n} \, \binom{4n}{2n}^2 \, }} \, ,
\label{eq:Lupas}
\end{equation}
which was implemented in \emph{Mathematica}$^\mathrm{TM}$ (version $6$) for computing $G$~\cite{Wolfram}.

On searching for new congruences modulo primes, Z.-W. Sun (2011) has conjectured that (see Conj.~A7 of Ref.~\cite{Sun})
\begin{equation}
\sum_{k=0}^{p-1}{(3 k +1) \, \frac{{\binom{2 k}{k}}^3}{(-8)^k}}  \overset{?} \equiv  p \left(\frac{-1}{p}\right) +p^3 \, E_{p-3} \, (\mathrm{mod} \,  p^4) \, ,
\label{eq:Sun0}
\end{equation}
which suggests that
\begin{equation}
\sum_{k=1}^{\infty}{(3 k -1) \: \frac{(-8)^k}{k^3 \: {\binom{2k}{k}}^3}}  \overset{?} =  -2 \, G  \, .
\label{eq:Sun}
\end{equation}
In fact, this result was confirmed computationally (with the aid of the WZ method) by Guillera in a recent work~\cite{Guillera}. The higher rate of convergence of this series in comparison to the previous ones, as well as its smaller number of basic arithmetic operations per iteration, has led me to search for a formal proof in order to justify its inclusion in further releases of mathematical softwares.

Here in this note, starting from a known hypergeometric identity I first prove a Ramanujan-type series representation for $G$ similar to that proved by Guillera. The convergence rate of this series is then shown to surpass that of all known similar central binomial series representations for this constant.  As a complement, I show that a simple proof of Eq.~\eqref{eq:Sun} follows from the application of a combinatorial identity to our main result.

\section{Ramanujan-type series for Catalan's constant}

Let us adopt the usual notation for the generalized hypergeometric series:
\begin{equation}
_{p}  F_{q} \! \left( \! \begin{array}{r} a_1 , \ldots , a_p \\ b_1 , \ldots, b_q \end{array} ; \, z \right) = \sum_{n=0}^\infty{\frac{\left(a_1\right)_n \, \ldots \, \left(a_p\right)_n}{\left(b_1\right)_n \, \ldots \, \left(b_q\right)_n} \, \frac{z^n}{n!}} \, ,
\end{equation}
where $\,\left(a\right)_n := a\,(a+1)\,\ldots\,(a+n-1) = \Gamma{(a+n)} / \Gamma{(a)}\,$ is the Pochhammer symbol for the rising factorial (by convention, $\left(a\right)_0 = 1$).  Our main result makes use of the lemma below, which determines a special value for $\, _{3}  F_{2} \! \left( \! \begin{array}{r} a_1 , a_2 , a_3 \\ b_1 , b_2 \end{array} ; \, z \right)$, a function that converges at $\,z=1\,$ whenever $\,\Re{\left\{(b_1+b_2)-(a_1+a_2+a_3)\right\}} > 0\,$ (see Eq.~(2.2.1) of Ref.~\cite{Slater}).  
\newline

\begin{lema}[A special value]
\label{lem:3F2}
\begin{equation*}
_{3}  F_{2} \! \left( \! \begin{array}{r} \frac12 , 1 , 1 \\ \frac32 , \frac32 \end{array} ; \, 1 \right) = 2\,G .
\label{eq:3F2}
\end{equation*}
\end{lema}

\begin{prova}
\; We start from a well-known integral representation for generalized hypergeometric functions (see, e.g., Eq.~(1.2) of Ref~\cite{Kratt}):
\begin{eqnarray}
_{p+1} F_{p} \! \left( \! \begin{array}{l} \alpha , \alpha_1 , \ldots , \alpha_p \\ \gamma , \beta_1 , \ldots , \beta_{p-1} \end{array} \! ; t \right) = \frac{\Gamma{(\gamma)}}{\,\Gamma{(\alpha)} \: \Gamma{(\gamma-\alpha)}} \nonumber \\
\times \int_0^1{z^{\alpha-1}(1-z)^{\gamma-\alpha-1} \, _{p}  F_{p-1} \! \left( \! \begin{array}{l} \alpha_1 , \ldots , \alpha_p \\ \beta_1 , \ldots , \beta_{p-1} \end{array} \! ;  t z \right) \, dz} \, ,
\end{eqnarray}
which holds whenever $\,\Re{(\alpha)}>0\,$ and $\,\Re{(\gamma-\alpha)}>0$. It then follows that 
\begin{eqnarray}
_{3}  F_{2} \! \left( \! \begin{array}{r} \frac12 , 1 , 1 \\ \frac32 , \frac32 \end{array} ; \, 1 \right) &=& \frac{\Gamma{(\frac32)}}{\Gamma{(\frac12)} \, \Gamma{(1)}} \, \int_0^1{z^{-\frac12}\,(1-z)^0 \, _{2}  F_{1} \! \left( \! \begin{array}{r} 1 , 1 \\ \frac32 \end{array} ; \, z \right) \, dz} \nonumber \\
&=& \frac12 \, \int_0^1{\frac{\,_{2}  F_{1} \! \left( \! \begin{array}{r} 1 , 1 \\ \frac32 \end{array} ; \, z \right)}{\sqrt{z}} \: dz} \nonumber \\
&=& \int_0^1{\,_{2}  F_{1} \! \left( \! \begin{array}{r} 1 , 1 \\ \frac32 \end{array} ; \, x^2 \right) \, dx} \, .
\label{eq:8}
\end{eqnarray}
Now, let us show that, for all non-null real $\,x \in (-1,1)$,
\begin{equation}
_{2} F_{1} \! \left( \! \begin{array}{r} 1 , 1 \\ \frac32 \end{array} ; \, x^2 \right) = \frac{\,\arcsin{x}}{\,x\,\sqrt{1-x^2}\:} \, .
\label{eq:pau}
\end{equation}
It is well-known that
\begin{equation}
_{2}  F_{1} \! \left( \! \begin{array}{r} 1/2 , 1/2 \\ \frac32 \end{array} ; \, x^2 \right) = \frac{\,\arcsin{x}}{x}
\label{eq:Slater}
\end{equation}
for all non-null values of $\,x\,$ for which the hypergeometric series at the left-hand side converges (see Eq.~(1.5.10) of Ref.~\cite{Slater}). Two successive applications of the Euler transformation formula, namely
\begin{equation}
_{2}  F_{1} \! \left( \! \begin{array}{r} a , b \\ c \end{array} ; \, z \right) = (1-z)^{-a} ~ _{2}  F_{1} \! \left( \! \begin{array}{r} a , c-b \\ c \end{array} ; \, \frac{z}{z-1} \right) ,
\end{equation}
on Eq.~\eqref{eq:Slater} lead us to
\begin{equation}
\frac{\arcsin{x}}{x} = \frac{1}{\sqrt{1-x^2}} ~ \, _{2}  F_{1} \! \left( \! \begin{array}{r} \frac12 , 1 \\ \frac32 \end{array} ; \, - \, \frac{x^2}{1-x^2} \right)
\end{equation}
and
\begin{equation}
\frac{\arcsin{x}}{x} = \frac{1-x^2}{\sqrt{1-x^2}} ~ \, _{2}  F_{1} \! \left( \! \begin{array}{r} 1 , 1 \\ \frac32 \end{array} ; \, x^2 \right) ,
\end{equation}
respectively. This completes the proof of Eq.~\eqref{eq:pau}.  From Eqs.~\eqref{eq:8} and \eqref{eq:pau}, one finds
\begin{equation}
_{3}  F_{2} \! \left( \! \begin{array}{r} \frac12 , 1 , 1 \\ \frac32 , \frac32 \end{array} ; \, 1 \right) = \int_0^1{\frac{\arcsin{x}}{x\,\sqrt{1-x^2}\:} \: d x} \, .
\end{equation}
The trigonometric substitution $\,x = \sin{\theta}\,$ reduces the last integral to
\begin{equation}
_{3}  F_{2} \! \left( \! \begin{array}{r} \frac12 , 1 , 1 \\ \frac32 , \frac32 \end{array} ; \, 1 \right) = \int_0^{\,\pi/2}{\!\frac{\theta}{\sin{\theta}} \: d\theta} \, .
\label{eq:trigsub}
\end{equation}
On noting that, apart from an arbitrary constant of integration,
\begin{equation}
\int{\frac{\theta}{\,\sin{\theta}} \: d\theta} = \theta \, \ln{\!\left[ -i \, \tan{\left(\frac{\theta}{2}\right)} \right]} +i \left[ \mathrm{Li}_2{\left(-e^{i \theta}\right)} - \mathrm{Li}_2{\left(e^{i \theta}\right)}\right] ,
\label{eq:primitiva}
\end{equation}
as can be easily checked by differentiating the function at the right-hand side, one finally finds
\begin{eqnarray}
\int_0^{\pi/2}{\frac{\theta}{\sin{\theta}} \, d\theta} = \frac{\pi}{2} \, \ln{\!\left( -i \right)} +i \left[ \, \mathrm{Li}_2{(-i)} - \mathrm{Li}_2{(i)}\right] \nonumber \\
-\left\{ \lim_{a \rightarrow 0^{+}}{a \, \ln{\!\left[ -i \, \tan{\left(\frac{a}{2}\right)} \right]}} +i \left[ \mathrm{Li}_2{(-1)} - \mathrm{Li}_2{(1)}\right]\right\} \nonumber \\
= \frac{\pi}{2} \, \left(\ln{1} -i\,\frac{\pi}{2}\right) +i \left[-2\,i \, \mathrm{Cl}_2{\left(\frac{\pi}{2}\right)}\right] -\left[0 +i \left( -\frac{\pi^2}{12} -\frac{\pi^2}{6}\right)\right] \nonumber \\
= -\,i\,\frac{\,\pi^2}{4} +2 \, \mathrm{Cl}_2{\left(\frac{\pi}{2}\right)} +i \, \frac{\,\pi^2}{4} \nonumber \\ = 2 \, G \, ,
\label{eq:2G}
\end{eqnarray}
where the special value of the Clausen function in Eq.~\eqref{eq:Clausen} and the principal value of the logarithm function, with $\,\mathrm{Arg}(z) \in (-\pi,\pi]$, were taken into account.\footnote{On Entry 9 of Adamchik's webpage~\cite{Adamchik}, where several representations for $\,G\,$ are proved computationally with \emph{Mathematica}$^\mathrm{TM}$, \, one finds $\, \frac12 \, \int_0^{\,\pi/2}{\frac{\theta}{\sin{\theta}} \, d\theta} = G$. Our Eqs.~\eqref{eq:primitiva} and \eqref{eq:2G} can then be viewed as a formal proof of this formula.}  The substitution of this result in Eq.~\eqref{eq:trigsub} completes the proof.
\begin{flushright} $\Box$ \end{flushright}
\end{prova}

We are now in a position to prove a rapidly converging central binomial series for the Catalan's constant, which is our main result.

\begin{teo}[Rapidly converging central binomial series]
\label{teo:main}
\begin{equation*}
G = \frac12 \, \sum_{n=0}^\infty{(-1)^n \, \frac{(3n+2) \: 8^n}{(2n+1)^3 \, {\binom{2n}{n}}^3}} \, .
\end{equation*}
\end{teo}

\begin{prova}
\; Firstly, note that the given series converges by the ratio test. Then, let
\begin{equation}
f(x) := \sum_{n=0}^\infty{(-1)^n \, \frac{\left(x+\frac12\right)_n^{\,3}}{8^n \, (x+1)_n^{\,3}} \, [6(x+n)+1]}
\end{equation}
be a function of a real variable $\,x\,$ in the open interval $\,(0,1)$. On noting that $\,(1)_n = n! \,$ and $\,\left(\frac32\right)_n = \dfrac{\Gamma\left(\frac32+n\right)}{\Gamma\left(\frac32\right)} = \dfrac{\left(n+\frac12\right) \, \Gamma\left(n+\frac12\right))}{\frac12 \, \Gamma(\frac12)} = \dfrac{(2n+1)!}{4^n \: n!}$, one finds
\begin{eqnarray}
f\!\left(\frac12\right) = \sum_{n=0}^\infty{(-1)^n \, \frac{(1)_n^{\,3} \: (6n+4)}{8^n \, \left(\frac32\right)_n^{\,3} }} = \sum_{n=0}^\infty{(-1)^n \, \frac{{n!}^3 \, (6n+4)}{8^n \, (2n+1)^3 \, \frac{{(2n)!}^3}{{n!}^3 \, 4^{3n}} } } \nonumber \\
= 2 \, \sum_{n=0}^\infty{(-1)^n \, \frac{{n!}^6 \, (3n+2) \, 64^n}{8^n \, (2n+1)^3 \, {(2n)!}^3 } } \nonumber  \\
= 2 \, \sum_{n=0}^\infty{(-1)^n \, \frac{(3n+2) \, 8^n}{(2n+1)^3 \, {\binom{2n}{n}}^3 } } \, .
\label{eq:fmeio}
\end{eqnarray}

On the other hand, from the third identity in Ref.~\cite{Guillera}, we know that
\begin{equation}
f(x) = 4 \,x \: \sum_{n=0}^\infty{ \frac{\left(x/2+\frac14\right)_n \, \left(x/2+\frac34\right)_n}{ {(x+1)_n}^2}} \,
\end{equation}
for all values of $\,x\,$ for which this series converges.  Therefore 
\begin{equation}
f(x) = 4 \, x \, ~ _{3}  F_{2} \! \left( \! \begin{array}{r} \frac{2x+1}{4} , \frac{2x+3}{4} , 1 \\ x+1 , x+1 \end{array} ; \, 1 \right) ,
\end{equation}
which implies that
\begin{equation}
f\!\left(\frac12\right) = 2 \; _{3}  F_{2} \! \left( \! \begin{array}{r} \frac12 , 1 , 1 \\ \frac32 , \frac32 \end{array} ; \, 1 \right) .
\end{equation}
From Lemma~\ref{lem:3F2}, this reduces to $\,f(\frac12) = 4\,G$.  From Eq.~\eqref{eq:fmeio}, we then find
\begin{equation}
2 \, \sum_{n=0}^\infty{(-1)^n \, \frac{(3n+2) \: 8^n}{(2n+1)^3 \, {\binom{2n}{n}}^3}} = 4 \, G \, .
\end{equation}

\begin{flushright} $\Box$ \end{flushright}
\end{prova}

As a bonus, we can use the above theorem to develop a formal proof for Eq.~\eqref{eq:Sun}.

\begin{corol}[Guillera's central binomial series]
\label{cor:Sun}
\begin{equation*}
G = -\frac12 \, \sum_{n=1}^\infty{(-1)^n \, \frac{(3n-1) \: 8^n}{n^3 \, {\binom{2n}{n}}^3}} \, .
\end{equation*}
\end{corol}

\begin{prova}
\; Since $\,\binom{2m+2}{m+1} = 2 \, \frac{2m+1}{m+1} \, \binom{2m}{m}\,$ for all non-negative integer values of $\,m$, then $\,(2m+1) \, \binom{2m}{m} = \frac12 \, (m+1) \, \binom{2m+2}{m+1}$. Theorem~\ref{teo:main} then yields
\begin{equation}
G = \frac12 \, \sum_{m=0}^\infty{(-1)^m \, \frac{(3m+2) \: 8^m}{\,\frac18\,(m+1)^3 \, {\binom{2m+2}{m+1}}^3 \,}} \, .
\end{equation}
The substitution $\:m = n -1\,$ completes the proof.

\begin{flushright} $\Box$ \end{flushright}
\end{prova}

The convergence rate of these novel central binomial series representations for $\,G\,$ is to be compared to that of other known similar series. This is done in details in the next section.

\section{Convergence rates}

Let us compare the convergence rates of the Ramanujan-type series for $\,G\,$ mentioned in this work.  By applying the Stirling's improved formula, namely
\begin{equation}
n! \sim \left(\frac{n}{e}\right)^n \sqrt{2 \pi \left(n+\frac16\right)} \, ,
\end{equation}
which implies that
\begin{equation}
\binom{2n}{n} \sim \frac{\,2^{2n}\sqrt{2 n+\frac16}}{\,\sqrt{2 \pi} \, \left( n+\frac16 \right)} \, ,
\end{equation}
we shall develop an order estimate of the $n$-th term for each series.
\newline

We begin with Ramanujan's series, given in Eq.~\eqref{eq:Rama}. Its $n$-th term is
\begin{eqnarray}
& & \frac{1}{(2n+1)^2 \, \binom{\,2n}{n}} \: \sim \: \frac{\sqrt{2 \pi \left(n+\frac16\right)}}{\, 2^{2n} \, \sqrt{2n+\frac16} \: (2n+1)^2} \nonumber \\
&=& \sqrt{\frac{\pi}{3}} \, \frac{6n+1}{\,2^{2n} \, \sqrt{12n+1} \: (2n+1)^2} 
\: \sim \: \frac{6n+1}{\, 2^{2n} \, \sqrt{12n+1} \: (2n+1)^2} \nonumber \\
& \sim & \frac{3}{\, 2^{2n} \, (2n+1) \, \sqrt{12n+1}} \, .
\end{eqnarray}

The more complex central binomial series by Lupas, see Eq.~\eqref{eq:Lupas}, has an $n$-th term whose absolute value can be estimated as follows:
\begin{eqnarray}
\frac{2^{8n} \, |40 n^2-24n +3|}{n^3 \, (2n-1) \, \binom{2n}{n} \, {\binom{4n}{2n}}^2} \, \sim \, \frac{2^{8n} \, |40 n^2-24n +3|}{n^3 \, (2n-1) \, \frac{2^{2n}\,\sqrt{2n+\frac16}}{\sqrt{2 \pi} \, (n+\frac16)} \, 2^{8n} \, \frac{3}{\pi} \, \frac{24n+1}{(12n+1)^2}} \nonumber \\
\sim \, \frac{\,|40 n^2-24n +3| \, (12n+1)^2 \, \sqrt{2 \pi} \, (n+\frac16)}{n^3 \, (2n-1) \, 2^{2n} \, \sqrt{2n+\frac16} \, (24n+1)} \, . 
\end{eqnarray}
For sufficiently large values of $\,n$, this simplifies to
\begin{eqnarray}
\frac{(40n-24) \, (12n+1)^2 \sqrt{2 \pi} \: (n+\frac16)}{\, 2^{2n} \, n^2 \, (2n-1) \, \sqrt{2n+\frac16} \: (24n+1)} 
\sim \frac{(40n-24) \, (12n+1)^2 \, (6n+1)}{\, 2^{2n} \, n^2 \, (2n-1) \, \sqrt{12n+1} \: (24n+1)} \nonumber \\
= \frac{4 \, (5n-3) \, (12n+1)^2 \, (6n+1)}{\, 2^{2n} \, n^2 \, (2n-1) \, \sqrt{12n+1} \: (12n+\frac12)} 
\sim \frac{4 \, (5n-3) \, (12n+1) \, (6n+1) \, \frac52}{\, 2^{2n} \, n^2 \, (5n-\frac52) \, \sqrt{12n+1}} \nonumber \\
\sim \frac{10 \, (12n+1) \, (6n+1)}{\, 2^{2n} \, n^2 \, \sqrt{12n+1}} 
\: = \: \frac{10 \, \sqrt{12n+1} \, (6n+1)}{\, 2^{2n} \, n^2} \nonumber \\
= \frac{5 \, \sqrt{12n+1} \, (6n+1)}{\, 2^{2n-1} \, n^2} \, .
\end{eqnarray}
This convergence is clearly slower than that of Ramanujan's series, not to mention that the number of basic arithmetic operations needed to compute the $\,n$-th term is considerably larger.\footnote{Despite these disadvantages, Lupas' series has been implemented in \emph{Mathematica}$^\mathrm{TM}$ (version $6$) for computing $\,G$~\cite{Wolfram}, which may be due to the absence of the term $\: \frac18 \, \pi \: \ln{\left(2+\sqrt{3}\,\right)}$, present in the Ramanujan's formula, which demands additional computational efforts.}

The $n$-th term of the central binomial series proved in our Theorem~\ref{teo:main} is
\begin{eqnarray}
\frac{(3n+2) \: 2^{3n}}{(2n+1)^3 \, {\binom{2n}{n}}^3} \, \sim \, \frac{(3n+2) \: 2^{3n}}{(2n+1)^3 \, \frac{2^{6n+3}}{(12n+2)^{\frac32}}} \nonumber \\
 =  \frac{3 (n+\frac23) \: 2^{3n} \, (12n+2)^{\frac32}}{8 \, (n+\frac12)^3 \, 2^{6n+3}} \, .
\end{eqnarray}
For sufficiently large values of $\,n$, this can be reduced to
\begin{eqnarray}
\frac38 \, \frac{(12n+2)^{\frac32}}{(n+\frac12)^2 \, 2^{3n+3}} \: = \: 3 \, \sqrt{2} \, \frac{(6n+1)^{\frac32}}{(2n+1)^2 \, 2^{3n+3}} \nonumber \\
\sim \: 9 \, \sqrt{2} \, \frac{\sqrt{6n+1}}{(2n+1) \, 2^{3n+3}} = 27 \, \sqrt{2} \, \frac{\sqrt{6n+1}}{\left(\sqrt{6n+3}\right)^2 \, 2^{3n+3}} \nonumber \\
\sim \: 27 \, \sqrt{2} \, \frac{1}{\sqrt{6n+3} \; 2^{3n+3}} \: \sim \: \frac{\,5 / {\sqrt{3}}}{\,2^{3n}  \, \sqrt{2n+1}} \: .
\end{eqnarray}
Note that the factor $\,2^{3n}\,$ in the denominator makes our series to converge faster than the Ramanujan's original series.
\newline

In Table~\ref{tabela}, below, we compare the error committed in approximating $\,G\,$ by the partial sum of the first $\,N\,$ terms corresponding to each central binomial series mentioned in this work. It is clear from this table that our series yields the smallest absolute error.

In a very recent paper by the Pilehroods~\cite{pile}, one finds in their Theorem~4 a rapidly converging central binomial series for $\,G$, namely
\begin{equation}
G = -\,\frac{1}{64} \, \sum_{n=1}^\infty{(-1)^n \, \frac{q(n) \; 256^n}{p(n) \: {\binom{8n}{4n}}^2 \: \binom{2n}{n} } }
\end{equation}
where $q(n)$ and $p(n)$ are polynomials in $n$ with degrees $6$ and $8$, respectively.  Although this series converges faster than the series in our Theorem~1, it was not included in our comparative study of Ramanujan-type series for $\,G\,$ because it presents polynomials with high degrees and large integer coefficients, as well as distinct binomial factors, which impedes us to characterize it as a Ramanujan-type series.  Moreover, it does not seem advantageous to use their series for high-precision computations of $\,G\,$ because it demands a greater number of basic operations per iteration in comparison to the simpler series proved here in this work.

Since a rational number cannot be approximated by rapidly converging sequences of (distinct) rational numbers and the Ap\'{e}ry's proof of irrationality for $\,\zeta(2)\,$ and $\,\zeta(3)\,$ is based upon the fast convergence of the central binomial series $\,\sum_{n=1}^\infty{1/\left[n^2\,\binom{2n}{n}\right]}\,$ and $\,\sum_{n=1}^\infty{(-1)^{n} / \left[n^3 \, \binom{2n}{n}\right]}$, respectively, as shown in Ref.~\cite{Poorten}, then it seems plausible that the series proved here could be taken as the starting point for a proof of irrationality for $G$.

\section*{Acknowledgments}
The author thanks Mrs.~Marcia R. Souza and Mr.~Bruno S. S. Lima for checking computationally all series in this work.  Thanks are also due to Dr.~Mathew Rogers for pointing out that a combinatorial identity could be taken into account to prove Eq.~\eqref{eq:Sun} directly from Theorem~\ref{teo:main}.
\newline

\newpage

\section*{Tables}

\begin{table}[ht]
\centering
\begin{tabular}{|l|l|l|l|}
\hline
$N$ & Lupas~\cite{Lupas} & Ramanujan~\cite{Ramanujan} &  Theorem~\ref{teo:main} \\ \hline
 &  &  & \\
10 & $2.0 \times 10^{-7}$ & $1.3 \times 10^{-9}$ & $3.1 \times 10^{-11}$ \\ 
 &  &  & \\
50 & $7.7 \times 10^{-32}$ & $1.1 \times 10^{-34}$ & $1.1 \times 10^{-47}$ \\ 
 &  &  & \\
100 & $4.3 \times 10^{-62}$ & $3.3 \times 10^{-65}$ & $5.6 \times 10^{-93}$ \\ 
 &  &  & \\
500 & $2.9 \times 10^{-303}$ & $4.6 \times 10^{-307}$ & $1.5 \times 10^{-454}$ \\ 
&  &  & \\
1000 & $1.9 \times 10^{-604}$ & $1.5 \times 10^{-608}$ & $2.9 \times 10^{-906}$ \\
 &  &  & \\ \hline
\end{tabular}
\caption{Absolute deviations from $\,G\,$ of the partial sums obtained by adding the first $\,N\,$ terms of each central binomial series mentioned in the text.}
\label{tabela}
\end{table}


\begin{thebibliography}{19}

\bibitem{Adamchik} V. S. Adamchik, \emph{33 representations for Catalan's constant}. Available at: \begin{verbatim} www.cs.cmu.edu/~adamchik/articles/catalan/catalan.htm \end{verbatim}


\bibitem{Bradley} D. M. Bradley, \emph{A class of acceleration formulae for Catalan's constant}, Ramanujan J. \textbf{3} (1999), pp.~159--173.

\bibitem{Catalan1865} E. Catalan, \emph{Sur la transformation des series, et sur quelques integrales definies}, Memoires de l'Academie Royale de Belgique (1865).

\bibitem{Finch} S. R. Finch, \textsl{Mathematical constants}, Cambridge University Press, Cambridge, 2003. Sec. 1.7.

\bibitem{Guillera} J. Guillera, \emph{Hypergeometric identities for 10 extended Ramanujan-type series}, Ramanujan J. \textbf{15} (2008), pp.~219--234.

\bibitem{pile} Kh. Hessami Pilehrood and T. Hessami Pilehrood, \emph{Series acceleration formulas for beta values}, Discrete Math. Theor. Comput. Sci. \textbf{12:2} (2010), pp.~223--236.

\bibitem{Kratt} C. Krattenthaler and K. Srinavasa Rao, \emph{Automatic generation of hypergeometric identities by the beta integral method}, J. Computat. Appl. Math. \textbf{160} (2003), 159--173.

\bibitem{Lima2012} F. M. S. Lima, \textsl{An Euler-type formula for $\,\beta(2n)\,$ and closed-form expressions for a class of zeta series}, Integral Transf. Spec. Funct. \textbf{23} (2012), 649--657. 
\bibitem{Lupas} A. Lupas, \emph{Formulae for some classical constants}. In \textsl{Proceedings of ROGER-2000} (2000). Available at:  \begin{verbatim} http://www.lacim.uqam.ca/~plouffe/articles/alupas1.pdf  \end{verbatim}

\bibitem{Poorten} A. van der Poorten, \emph{A proof that Euler missed: Ap\'{e}ry's proof of the irrationality of $\,\zeta(3)$}, Math. Intelligencer \textbf{1} (1979), pp.~195--203.

\bibitem{Ramanujan} S. Ramanujan, \emph{On the integral $\int_0^x{\frac{\tan^{-1}{t}}{t} \, dt}$}, J. Indian Math. Soc. \textbf{VII} (1915), pp.~93--96.

\bibitem{RZ} T. Rivoal and W. Zudilin, \emph{Diophantine properties of numbers related to Catalan's constant}, Math. Annalen \textbf{326} (2003), pp.~705--721.

\bibitem{Slater} L. J. Slater, \textsl{Generalized Hypergeometric Functions}, Cambridge University Press, Cambridge, 1966. 

\bibitem{Sun} Z.-W. Sun, \emph{Open conjectures on congruences}. arXiv:math.NT/0911.5665v59 (2011).

\bibitem{Wolfram} O. Marichev, J. Sondow, and E. W. Weisstein, \emph{Catalan's constant}. From MathWorld, a Wolfram web resource. Available at:  \begin{verbatim} http://mathworld.wolfram.com/CatalansConstant.html  \end{verbatim}

\bibitem{Yee}  A. Yee, \emph{Large Computations}, 20 June 2012. Available at:  \begin{verbatim} http://www.numberworld.org/digits/Catalan
\end{verbatim}

\end{thebibliography}
\end{document}